\documentclass[a4paper, 12pt]{article}
\usepackage{adjustbox}
\usepackage{graphicx}
\usepackage{amsfonts}
\usepackage{amsbsy}
\usepackage{amssymb}
\usepackage[fleqn]{amsmath}
\usepackage{multicol}
\usepackage{longtable,ltxtable,booktabs}
\usepackage{blkarray}
\usepackage{afterpage}
\usepackage{float}
\usepackage{caption}
\usepackage{subcaption}
\usepackage{multirow}
\usepackage{pdflscape}
\usepackage{latexsym}
\usepackage{rotating}
\usepackage{enumerate}

\usepackage{mathtools}

\usepackage{epstopdf}
\usepackage{caption}
\usepackage{xcolor}
\usepackage{amsmath}
\usepackage{longtable}
\usepackage[margin=.8 in]{geometry}

\usepackage{verbatim}

\usepackage{mathrsfs}
\usepackage{amsthm}
\usepackage{setspace}

 \usepackage{pdflscape}
\usepackage{colortbl}
\usepackage{latexsym}
\usepackage{amsfonts}
\usepackage{epsfig}
\usepackage{hyperref}
\usepackage{graphics}
\usepackage{float}
\theoremstyle{definition}
\newtheorem{thm}{Theorem}[section]

\newtheorem{defn}{Definition}[section]
\newtheorem{prop}{Proposition}[section]

\newtheorem{lemma}{Lemma}[section]

\newtheorem{cor}{Corollary}[section]

\theoremstyle{definition}

\title{p-Groups in which kernels of the non-linear irreducible characters are of equal order }

\date{{}}\author{Nabajit Talukdar\thanks{Corresponding author.Email address : ntalukdar2000@yahoo.co.in}\\
	\small Department of Mathematics\\
	\small Cotton University\\
	\small Guwahati-781001, India}

\begin{document}
\maketitle

\begin{abstract}For an irreducible character $\chi$ of a finite group $G$, its kernel is defined as $\text{ker }\chi=\{g\in G: \chi(g)=\chi(1)\}$. 
In this paper we characterize the finite groups of prime power order(for odd prime) in which kernels of all of the non-linear irreducible characters are of the same order.
\end{abstract}

  $\mathbf{2010\;Mathematics\;Subject\;Classification:}\;20C15$
	
	$\mathbf{Keywords\;and\;Phrases:} $ p-group, character kernel, generalized Camina pair

\section{Introduction}
In this paper, all groups  are finite. By a $p$-group we denote a group having order some power of $p$, where $p$ is an odd prime. By $\text{Irr}(G)$ and $\text{Irr}_{1}(G)$  we denote the set of of  complex irreducible characters and the set of   complex non-linear irreducible characters of the group $G$ respectively. By $c(G)$ we denote nilpotency class of $G$. For the terminologies not defined here follow the Isaacs' book \cite{isaacs1994character} and Rotman's book \cite{rotman2012introduction}.

\begin{defn} For any finite group $G$, we define:
\begin{enumerate}[(i)]
\item $\text{Kern}(G)=\{\text{ker} \chi: \chi \in \text{Irr}_{1}(G)\}$.
\item $\text{sk}(G)=\{|\text{ker} \chi|: \chi \in \text{Irr}(G)\}$.
\item $\text{skn}(G)=\{|\text{ker} \chi|: \chi \in \text{Irr}_{1}(G)\}$.
\end{enumerate}
\end{defn}

 There are many research articles(e.g. \cite{wang2008note}, \cite{doostie2012finite}, \cite{li2023finite}) where studies were done on kernels of irreducible characters of a group. In this paper we we characterize the finite groups of prime power order(for odd prime) in which kernels of all of the non-linear irreducible characters are of the same order. Following \cite{nenciu2012isomorphic}, 
a non-Abelian group $G$ is called a generalized VZ-group (GVZ for short) if for any $\chi\in \text{Irr}(G)$ we have that $\chi(g)=0$ for all $g\in G\setminus Z(\chi)$. By [Corollary 2.30 \cite{isaacs1994character}] we get that if $G$ is a GVZ-group then $\chi(1)^{2}=|G/Z(\chi)|$. A group of nilpotency class $2$ is obviously a GVZ-group. Studies of GVZ-groups having two character degrees were done by K. K. Rajkhowa and N. Talukdar{\cite{talukdar2025gvz}. It turns out that a finite group $G$ of prime power order(for odd prime) in which kernels of all of the non-linear irreducible characters are of the same order is a GVZ-group having two character degrees if the $|G'|$ does not exceed the order of the kernel of any of the non-linear irreducible characters of $G$.

The main result of this paper is the following.

  \begin{thm}
\label{thm_kernel_same}
     Let $G$ be a non-Abelian $p$-group, where $p$ is an odd prime. Then kernel of each of the non-linear irreducible characters of $G$ has order $p^{m}$, $m\geq 1$ if and only if  $c(G)=2$ and one of the following holds:
\begin{enumerate}[1.]
      \item $G'\cong C_{p}\times C_{p} \times \ldots \times C_{p}(\ (m+1) \text{ times})$ is the unique normal subgroup of $G$ of order $p^{m+1}$.
     \item  $\text{cd}(G)=\{1, (\frac{|G|}{p^{m+1}})^{\frac{1}{2}}\}$ and $Z(G/(G'\cap \text{ker }\chi))$ is elementary Abelian for every $\chi\in \text{Irr}_{1}(G)$ ($|G'|\leq p^{m}$ in this case). 
     \end{enumerate}
  \end{thm}

\section{Preliminaries}
In this section we shall prove some preliminary results which will aid us in proving our main theorem. First we state a result that says that any normal subgroup of a finite group can be obtained by taking the intersection of the kernels of some of its irreducible characters. 
\begin{lemma}
\label{lemma_normal_intersect_kernel}
\cite{alperin1995group}
Let $K_{1}, K_{2}, \ldots , K_{r}$ be the kernels of the irreducible characters of a group $G$. If $N\unlhd G$, then $N=\cap_{i\in I}K_{i}$ where $I\subseteq \{1,2, \ldots , r\}$.
\end{lemma}

The following result says that the intersection of the kernels  of the non-linear irreducible characters of a group is trivial.
\begin{lemma}
\label{lemma_non_lin_inter}
\cite{li2019finite}
For a group $G$, $\cap_{\chi \in \text{Irr}_{1}(G)}\text{ker}\ \chi=\{1\}$.
\end{lemma}

\begin{lemma}
\label{lemma_centre_kernel}
    Let $G$ be a group and $\chi \in \text{Irr}(G)$. Then $[Z(\chi), G]\leq \text{ker}\ \chi$.
    \begin{proof}
        This holds since $Z(\chi)/\text{ker}\ \chi=Z(G/\text{ker}\ \chi)$.
    \end{proof}
\end{lemma}

In the following Lemma we state the necessary and sufficient condition that all the non-linear irreducible characters of a $p$-group are faithful.
 \begin{lemma}
 \label{lemma_all_faithful}
  \cite{doostie2012finite}
  Let $G$ be a $p$-group. Then $\text{skn}(G)=\{1\}$ if and only if $|G'|=p$ and $Z(G)$ is cyclic.
 \end{lemma}

In the following result we obtain an upper bound for the nilpotency class of a $p$-group.
\begin{lemma}
Let $G$ be a non-Abelian $p$-group and $m$ be the largest positive integer such that  $p^{m} \in \text{skn}(G)$. Then $c(G)\leq m+2$.
\begin{proof}
   Since $m$ be the largest positive integer such that  $p^{m} \in \text{skn}(G)$, it follows that if $N\unlhd G$ and $|N|\geq p^{m+1}$, then $G'\leq N$. If $c(G)\leq m$, we are done. Suppose $c(G)> m$. Then it follows that $|Z_{m+1}(G)|\geq p^{m+1}$ and hence $G'\leq Z_{m+1}(G)$. Thus we get that $c(G)\leq m+2$.
\end{proof}
\end{lemma}

\begin{defn}
    \cite{fernandez2001groups}
    We say that a $p$-group $G$ satisfies the strong condition on normal
subgroups provided that, for any $N\unlhd G$, either  $G'\leq N$ or $N\leq Z(G)$. Again, we say that a $p$-group $G$ satisfies the weak condition on normal subgroups provided that, for any $N\unlhd G$, either  $G'\leq N$ or $|NZ(G)/Z(G)|\leq p$.
 
\end{defn}

\begin{thm}
   \cite{fernandez2001groups}
   \label{thm_nc_2_str_weak}
   Let $G$ be a $p$-group of nilpotency class $2$.
   \begin{enumerate}[(i)]
       \item If $G$ satisfies the strong condition on normal subgroups, then $\text{exp }G/Z(G)=\text{exp } G'=p$.
       \item If $G$ satisfies the weak condition on normal subgroups, then $\text{exp }G/Z(G)=\text{exp } G'=p$ or $p^{2}$. In the latter case, $G/Z(G)\cong C_{p^{2}}\times C_{p^{2}}$ and $G'\cong C_{p^{2}}$.
   \end{enumerate}
\end{thm}

Now we state a result regarding a $p$-group in which $|G'|=p$.
\begin{lemma}
\label{lemma_derived_p}
\cite{talukdar2025small}
Let $G$ be a $p$-group such that $|G'|=p$. Then 
\begin{enumerate}[1.]
\item cd$(G)=\{1, |G/Z(G)|^{\frac{1}{2}}\}$.
\item $G$ satisfies the strong condition on normal subgroups.
\end{enumerate}
\end{lemma}

Following \cite{zhmud1999characters}, two subgroups $F$ and $H$ of a a group $G$ are said to be nonincident if $F\neq F\cap H\neq H$. A non-Abelian group is said to be a $J$-group if any two different elements of $\text{Kern}(G)$ are nonincident. For a $p$-group $G$, if $\text{skn}(G)=\{p^{m}\}$, $m\geq 1$, then the group $G$ is a $J$-group. For a group $G$, let $\mathcal{N}(G)=\{N \ |\ N\unlhd G, G'\not \subseteq N\}$. Following \cite{li2023finite}, a group $G$ is said to be a $K(m,n)$ group if $|\text{Kern}(G)|=m$ and $|\mathcal{N}(G)|=n$.

\begin{lemma}
\cite{zhmud1999characters}
    A $p$-group is a $J$-group if and only if $G'\subseteq Z(G)$ and $\text{exp}\ G'=p$.
\end{lemma}

\begin{cor}
    Let $G$ be a $p$-group such that $\text{skn}(G)=\{p^{m}\}$, $m\geq 1$. Then $\text{exp}\ G'=p$.
\end{cor}

In the following result we obtained the orders of the kernels of the irreducible characters of an Abelian $p$-group.

\begin{lemma}
\label{lemma_ker_size_abelian}
Let $G$ be an Abelian $p$-group such that $|G|=p^{n}$ and $\text{exp}(G)=p^{e}$. Then $sk(G)=\{p^{n-e}, \ldots, p^{n}\}$. In particular if $G$ is elementary Abelian, then $sk(G)=\{p^{n-1}, p^{n}\}$.
\begin{proof}
   From [ Problem 2.7, \cite{isaacs1994character}], we get that $\text{Irr}(G)$ is isomorphic to $G$. Thus we get that $\{|G/ker\ \chi| : \chi \in \text{Irr}(G)\}=\{1, p, \ldots , p^{e}\}$. Hence the result follows.
\end{proof}
\end{lemma}

Now we prove some results in case of a non-Abelian $p$-group.

\begin{lemma}
\label{lemma_derived_subg_order}
Let $G$ be a non-Abelian $p$-group. If if  $m$ be the largest positive integer such that  $p^{m} \in \text{skn}(G)$, then either $G'$ is the unique normal subgroup of order $p^{m+1}$ or $|G'|\leq  p^{m}$.
\begin{proof}  
    Let $K$ be the kernel of a non-linear irreducible character of $G$ such that $|K|=p^{m}$. Now all the non-linear irreducible characters of the group  $G/K$ are faithful. Hence by Lemma \ref{lemma_all_faithful} $|G'/(G'\cap K)|=p$. This gives that $|G'|\leq p|K|=p^{m+1}$. Suppose $|G'|=p^{m+1}$. Let $N\unlhd G$ be such that $|N|=p^{m+1}$. Since $G/N$ is Abelian we get that $G'\subseteq N$ and hence  $N=G'$. This proves that $G'$ is the unique normal subgroup of order $p^{m+1}$.   
\end{proof}
\end{lemma}

\begin{lemma}
  \label{lemma_degree_square}
    Let $G$ be a non-Abelian $p$-group and $m$ be the largest positive integer such that  $p^{m} \in \text{skn}(G)$. Let $\chi \in \text{Irr}_{1}(G)$ be such that $|\text{ker}\ \chi|=p^{m}$. Then $\chi(1)^{2}=|G/Z(\chi)|$.
    \begin{proof}
        Since $|Z(\chi)|>p^{m}$, it follows from Lemma \ref{lemma_normal_intersect_kernel} that $G'\subseteq Z(\chi)$. From [Theorem 2.31, \cite{isaacs1994character}] it follows that $\chi(1)^{2}=|G/Z(\chi)|$.
    \end{proof}
\end{lemma}

\begin{lemma}
\label{lemma_cent_bigger_ker}
  Let $G$ be a $p$-group and $\chi \in \text{Irr}(G)$ be such that $\chi \neq 1_{G}$. Then $\text{ker}\ \chi<Z(\chi)$.
  \begin{proof}
      We suppose that $\text{ker}\ \chi=Z(\chi)$. Then $Z(G/\text{ker}\ \chi)=Z(\chi)/\text{ker}\ \chi=1$. This is possible only if the group $G/\text{ker}\ \chi$ is the trivial group. Thus $\text{ker}\ \chi=G$ and hence it follows that $\chi = 1_{G}$. Thus if $\chi \neq 1_{G}$, $\text{ker}\ \chi<Z(\chi)$. 
  \end{proof}
\end{lemma}

\begin{thm}
\label{thm_kern_1_p}
    For a $p$-group $G$, $\text{skn}(G)=\{1,p\}$ if and only if $G$ is a group of order $p^{4}$ of maximal class.

\end{thm}

The structures of Abelian groups have been discussed in \cite{kurzweil2004finite} and \cite{rotman2012introduction}. In particular we state the following results.

\begin{prop}
\label{prop_kurz_abelian}
\cite{kurzweil2004finite}
Let $G$ be a finite Abelian group and $U$ be a cyclic subgroup of maximal order in $G$. Then there exists a complement $V$ of $U$ in $G$.
\end{prop}

\begin{lemma}
\label{lemma_abelian_quotient_subgroup}
\cite{rotman2012introduction}
Let $G$ be a finite Abelian group and $H\leq G$. Then $G$ contains a subgroup isopmorphic to $G/H$.
\end{lemma}

A pair $(G,N)$ is said to be a generalized Camina pair (abbreviated GCP) if $N$ is a normal subgroup of the group $G$ and all the non-linear irreducible characters of $G$ vanish outside $N$. 
The notion of GCP was introduced by Lewis in \cite{lewis2009vanishing}. An equivalent condition for a pair $(G,N)$ to be a GCP is: A pair $(G,N)$ is a GCP if and only if for $g\in G\setminus N$, the conjugacy class of $g$ in $G$ is $gG'$.
The following theorem can be obtained from the statement and proof of [Theorem 3.1, \cite{prajapati2017irreducible}].

\begin{thm}
\label{thm_gcp_praja}
 
Let $(G,Z(G))$ be a GCP. Then we have the following.

   \begin{enumerate}[(i)]

       \item $cd(G)=\{1, |G/Z(G)|^{\frac{1}{2}}\}$.

       \item There is a bijection between the sets $\text{Irr}(Z(G)|G')$ and $\text{Irr}_{1}(G)$ and the bijection is given by $\lambda \rightarrow \chi$, where
       \[
\chi(g)=
\begin{cases}
0 & \text{if} \ g \not \in Z(G),\\
|G/Z(G)|^{\frac{1}{2}}\lambda(g) & \text{if} \ g \in Z(G).\\
  
\end{cases}
\]
       
   \end{enumerate}
\end{thm}

\begin{lemma}
\label{lemma_g'_gcp}
    Let $G$ be a $p$-group such that $|G'|=p$. Then $(G, Z(G))$ is a generalized Camina pair.
    \begin{proof}
        We shall show that $Z(\chi)=Z(G)$ for all $\chi \in \text{Irr}_{1}(G)$.
If possible let $\chi \in \text{Irr}_{1}(G)$ be such that $Z(\chi)\neq Z(G)$. We choose $g\in Z(\chi) \setminus Z(G)$.
Then there exists $ h\in G$ such that $gh\neq hg$ that is $[g,h](\neq 1)\in G'$.
Since $|G'|=p$, it follows that $G'$ is cyclic and $G'=<[g,h]>$. Let $\Psi$ be the representation of $G$ that affords the character $\chi$. Since $g\in Z(\chi)$, $\Psi(g)=\epsilon I$ for some $\epsilon \in \mathbb{C}$. Thus we have that $\Psi(g)\Psi(h)=\Psi(h)\Psi(g)$ and hence  $[g,h]\in \text{ker}\ \Psi=\text{ker}\ \chi$. Then it follows that $G'=<[g,h]> \leq ker\chi$. Consequently we get that $\chi$ is a linear character. This contradiction proves that $Z(G)=Z(\chi)$.
Since $G'\subseteq Z(G)=Z(\chi)$, $G/Z(\chi)$ is Abelian. By[Corollary  2.30 and Theorem 2.31, \cite{isaacs1994character}] we get that $(G, Z(G))$ is a generalized Camina pair.
    \end{proof}
\end{lemma}

\section{Main Results}
We prove the Theorem \ref{thm_kernel_same} by means of the following Lemmas and Propositions.
\begin{prop}
\label{prop_size_1_all}
     Let $G$ be a non-Abelian $p$-group and let $\text{skn}(G)=\{p^{m}\}$, $m\geq 1$. Then 
     \begin{enumerate}[1.]
       \item $c(G)=2$.
       \item $G/Z(G)$ is elementary-Abelian.
     \end{enumerate}

     \begin{proof}
       \begin{enumerate}[1.]
       
       \item From Lemma \ref{lemma_cent_bigger_ker} it follows that $G'\subseteq Z(\chi)$ for all $\chi \in \text{Irr}_{1}(G)$. Hence $G'\subseteq \cap_{\chi \in \text{Irr}_{1}(G)}Z(\chi)=Z(G)$. Thus 
       $c(G)=2$.

       \item For any $\chi \in \text{Irr}_{1}(G)$ let $K=\text{ker }\chi$. We get that $G/K$ is a group such that every non-linear irreducible character of $G/K$ is faithful. Thus $(G/K)'$ is the unique minimal normal subgroup of $G/K$. By [Lemma 12.3, \cite{isaacs1994character}] we get that $G/Z(\chi)$ is elementary Abelian. Hence for any $x\in G$, $x^{p}\in Z(\chi)$. Thus we get that $x^{p}\in \cap_{\chi \in \text{Irr}_{1}(G)}Z(\chi)=Z(G)$. This gives that $G/Z(G)$ is elementary-Abelian.
      \end{enumerate}
     \end{proof}
\end{prop}

\begin{lemma}
\label{lemma_g'_unique-eqker}
   Let $G$ be a non-Abelian $p$-group. Let $\text{skn}(G)=\{p^{m}\}$, $m\geq 1$ and $G'$ is the unique normal subgroup of $G$ of order $p^{m+1}$. Then $G'$   is elementary Abelian.
   \begin{proof}
        We get that $K\subseteq G' \subseteq Z(G)$ for all $K\in \text{Kern}(G)$. Hence $G$ satisfies the strong condition on normal subgroups. Thus from Theorem \ref{thm_nc_2_str_weak} we get that $\text{exp}\ G'=p$. Hence $G'$   is elementary Abelian.
   \end{proof}
\end{lemma}

\begin{lemma}
\label{lemma_g'_p_kernel_size}
    Let $G$ be a $p$-group with $|G'|=p$. If $\text{skn}(G)=\{p^{m}\}$,  $m\geq 0$, then $\frac{|Z(G)|}{\text{exp}\ Z(G)}=p^{m}$.
        \begin{proof}
       We prove the result by induction on $m$. If $m=0$, each of the non-linear irreducible characters of $G$ is faithful and hence $Z(G)$ is cyclic. Thus we get that $\frac{|Z(G)|}{\text{exp}\ Z(G)}=\frac{|Z(G)|}{|Z(G)|}=p^{0}$. Since $|G'|=p$, $G$ satisfies the strong condition on normal subgroups and hence $K\leq Z(G)$ for all $K\in \text{Kern}(G)$. 
        Now we consider that $m\geq 1$. From Lemma \ref{lemma_all_faithful} we get that $Z(G)$ is not cyclic. First we assume that $\text{skn}(G)=\{p\}$. By Lemma \ref{lemma_abelian_quotient_subgroup} we get that for any $K\in \text{Kern}(G)$, $Z(G)$ contains a subgroup isomorphic to $Z(G)/K$. Since $Z(G)/K$ is cyclic and $|Z(G)/K|=\frac{|Z(G)|}{p}$, we get that $\text{exp}\ Z(G)=\frac{|Z(G)|}{p}$. This proves that $\frac{|Z(G)|}{\text{exp}(Z(G))}=p$. Now we assume that the result holds if $\text{skn}(G)=\{p^{m-1}\}$.
        Suppose $p^{e}=\text{exp}\  Z(G)$ and $a\in Z(G)$ be such that $o(a)=p^{e}$. If $a^{p^{e-1}}\in K$ for all $K\in \text{Kern}(G)$, then by Lemma \ref{lemma_non_lin_inter}, we get that $a^{p^{e-1}}=1$. This contradicts that $o(a)=p^{e}$. Thus we can choose  $K\in \text{Kern}(G)$ such that $a^{p^{e-1}}\not \in K$. Let $N\subseteq K\cap Z(G)$ be of order $p$. Then $\text{skn}(G/N)=\{p^{m-1}\}$. By induction hypothesis we get that $\frac{|Z(G/N)|}{\text{exp}\ Z(G/N)}=p^{m-1}$. From Lemma \ref{lemma_derived_p} and [\cite{fernandez2001groups}, Theorem A] we get that $Z(G/N)=Z(G)/N$. Since $a^{p^{e-1}}\not \in N$, we get that $p^{e}=\text{exp}\  Z(G)/N$. Thus we get that $\frac{|Z(G)/N|}{\text{exp}\ Z(G)}=p^{m-1}$ and hence $\frac{|Z(G)|}{\text{exp}\ Z(G)}=p^{m}$.

    \end{proof}
\end{lemma}

\begin{lemma}
\label{lemma_g'_p_centre}
    Let $G$ be a $p$-group with $|G'|=p$. Then $\text{skn}(G)=\{p\}$ if and only if $Z(G)\cong C_{p}\times C_{p}$.
    \begin{proof}
        We get that $\text{exp}\ Z(G)=\frac{|Z(G)|}{p}$. By Proposition \ref{prop_kurz_abelian} we get that $Z(G)\cong C_{p^{k-1}}\times C_{p}$, where $|Z(G)|=p^{k}$. Thus $Z(G)$ contains $p+1$ subgroups of order $p$, one of which is $G'$ and the remaining  are the members of $\text{Kern}(G)$. Hence $G$ is $K(p,p+1)$ group. From [Theorem 4.2. \cite{li2023finite}] we get that  $Z(G)\cong C_{p}\times C_{p}$.\\
        Conversely, suppose that $Z(G)\cong C_{p}\times C_{p}$. Then $Z(G)$ contains $p+1$ subgroups of order $p$, one of which is $G'$. Since $|G'|=p$, $G$ satisfies the strong condition on normal subgroups. Thus $K\leq Z(G)$ for every $K\in \text{Kern}(G)$. This gives that $\text{skn}(G)=\{p\}$.
    \end{proof}
    \end{lemma}

\begin{prop}
\label{prop_skng_pm}
    For a $p$-group $G$, $\text{skn}(G)=\{p\}$ if and only if $c(G)=2$ and one of the following holds:
    \begin{enumerate}[1.]
     \item $|G'|=p$ and $Z(G)\cong C_{p}\times C_{p}$.

     \item $G'\cong C_{p}\times C_{p}$ is the unique normal subgroup of order $p^{2}$ of $G$.
    \end{enumerate}    
    \begin{proof}
      Let $\text{skn}(G)=\{p\}$. From Proposition \ref{prop_size_1_all} we get that  $c(G)=2$. By Lemma \ref{lemma_derived_subg_order} we get that $|G'|\leq p^{2}$. If $|G'|=p$, from Lemma \ref{lemma_g'_p_centre} we get that $Z(G)\cong C_{p}\times C_{p}$. If $|G'|=p^{2}$, by Lemma \ref{lemma_derived_subg_order} we get that $G'$ is the unique normal subgroup of order $p^{2}$. Since $\text{skn}(G)=\{p\}$, $G$ satisfies the strong condition on normal subgroups. Thus from Theorem \ref{thm_nc_2_str_weak} we get that $\text{exp}\ G'=p$ and consequently $G'\cong C_{p}\times C_{p}$. \\
      Conversely suppose that $c(G)=2$ and either $|G'|=p$ and $Z(G)\cong C_{p}\times C_{p}$ or $G'\cong C_{p}\times C_{p}$ is the unique normal subgroup of order $p^{2}$ of $G$. If $|G'|=p$ and $Z(G)\cong C_{p}\times C_{p}$ from Lemma \ref{lemma_g'_p_centre} we get that $\text{skn}(G)=\{p\}$. Now let $G'\cong C_{p}\times C_{p}$ is the unique normal subgroup of order $p^{2}$ of $G$. Let $K$  be the kernel of a non-linear irreducible character of $G$. Then $K <G'$ and hence $|K|\in \{1,p\}$. Thus $\text{skn}(G)=\{1\}$ or $\text{skn}(G)=\{p\}$ or $\text{skn}(G)=\{1, p\}$. From Lemma \ref{lemma_all_faithful} and Theorem \ref{thm_kern_1_p}, we can rule out the options $\text{skn}(G)=\{1\}$ and $\text{skn}(G)=\{1, p\}$. Hence, it follows that $\text{skn}(G)=\{p\}$. 
    \end{proof}
    \end{prop}

    \begin{prop}
    \label{prop_g'_p_centre}
    Let $G$ be a $p$-group such that $|G'|=p$. Then $\text{skn}(G)=\{p^{m}\}$, $m\geq 1$ if and only if $Z(G)\cong C_{p}\times C_{p}\times \ldots \times C_{p}(\ (m+1) \text{ times})$.
    \begin{proof}
      Let  $|G'|=p$ and $\text{skn}(G)=\{p^{m}\}$. If $\text{skn}(G)=\{p\}$, we get that $Z(G)\cong C_{p}\times C_{p}$. Next let $\text{skn}(G)=\{p^{m}\}$ and   $m\geq 2$. Let $K\in \text{Kern}(G)$ and let $K/N$ be a chief factor of $G$. Then $|(G/N)'|=p$ and $\text{skn}(G/N)=\{p\}$. Hence $|Z(G/N)|=p^{2}$. By Lemma \ref{lemma_derived_p}
and by [Theorem A, \cite{fernandez2001groups}] we get that $Z(G/N)=Z(G)/N$. Thus $|Z(G)/N|=p^{2}$ and hence $|Z(G)/K|=p$. By Lemma \ref{lemma_g'_p_kernel_size} we get that $\text{exp}\ Z(G)=\frac{|Z(G)|}{|K|}=p$. Thus it follows that $Z(G)\cong C_{p}\times C_{p}\times \ldots \times C_{p}(\ (m+1) \text{ times})$.
        
        Conversely we assume that $|G'|=p$ and $Z(G)\cong C_{p}\times C_{p}\times \ldots \times C_{p}(\ (m+1) \text{ times})$. From Lemma \ref{lemma_ker_size_abelian} we get that the kernel of every member of $\text{Irr}(Z(G)|G')$ has order $p^{m}$. Now $(G, Z(G))$ is a generalized Camina pair. From Theorem \ref{thm_gcp_praja} we get that $\text{skn}(G)=\{p^{m}\}$.
    \end{proof}
    \end{prop}

In the following theorem we prove our main results.

  \begin{thm}
\label{thm_derived_unique_0}
   Let $G$ be a non-Abelian $p$-group, where $p$ is an odd prime. Then kernel of each of the non-linear irreducible characters of $G$ has order $p^{m}$, $m\geq 1$ if and only if  $c(G)=2$ and one of the following holds:
\begin{enumerate}[1.]
      \item $G'\cong C_{p}\times C_{p} \times \ldots \times C_{p}(\ (m+1) \text{ times})$ is the unique normal subgroup of $G$ of order $p^{m+1}$.
     \item  $\text{cd}(G)=\{1, (\frac{|G|}{p^{m+1}})^{\frac{1}{2}}\}$ and $Z(G/(G'\cap \text{ker }\chi))$ is elementary Abelian for every $\chi\in \text{Irr}_{1}(G)$ ($|G'|\leq p^{m}$ in this case). 
 \end{enumerate}

 \begin{proof}
 Let $\text{skn}(G)=\{p^{m}\}$, $m\geq 1$. From Proposition \ref{prop_size_1_all} we get that $c(G)=2$. By Lemma \ref{lemma_derived_subg_order}, we get that either $G'$ is the unique normal subgroup of order $p^{m+1}$ or $|G'|\leq  p^{m}$. If $G'$ is the unique normal subgroup of order $p^{m+1}$, from Lemma \ref{lemma_g'_unique-eqker} we get that $G'\cong C_{p}\times C_{p} \times \ldots \times C_{p}$. Next we consider $|G'|\leq  p^{m}$.  Let $\chi \in \text{Irr}_{1}(G)$. If $G'\cap \text{ker }\chi=\text{ker }\chi$, then $\text{ker }\chi< G'$ and hence $|G'|\geq p^{m+1}$. Thus we have that  $G'\cap \text{ker }\chi<\text{ker }\chi$ and consequently $|\text{skn}(G/(G'\cap \text{ker }\chi))|=1$. Since $|(G/(G'\cap \text{ker }\chi))'|=p$, by Proposition \ref{prop_g'_p_centre} we get that $Z(G/(G'\cap \text{ker }\chi))$ is elementary Abelian. Since $[Z(\chi), G]\subseteq G'\cap \text{ker }\chi$, we get that $Z(\chi)/(G'\cap \text{ker }\chi)\subseteq Z(G/(G'\cap \text{ker }\chi))$. Consequently $Z(\chi)/(G'\cap \text{ker }\chi)$ is elementary Abelian and hence $Z(\chi)/ \text{ker }\chi$ is elementary Abelian. Since $Z(\chi)/ \text{ker }\chi$ is cyclic and $ \text{ker }\chi< Z(\chi)$, we get that $Z(\chi)/ \text{ker }\chi\cong C_{p}$. Thus $|Z(\chi)|= p^{m+1}$ and consequently  from Proposition \ref{prop_size_1_all} we get that $\chi(1)^{2}=|G/Z(\chi)|$ and hence $\text{cd}(G)=\{1, (\frac{|G|}{p^{m+1}})^{\frac{1}{2}}\}$.\\
For the converse part, first we assume that  $c(G)=2$  and $G'\cong C_{p}\times C_{p} \times \ldots \times C_{p}$ is the unique normal subgroup of $G$ of order $p^{m+1}$. Then it is easy to see that $\text{skn}(G)=\{p^{m}\}$. Next we assume that $c(G)=2$, $\text{cd}(G)=\{1, (\frac{|G|}{p^{m+1}})^{\frac{1}{2}}\}$ and $Z(G/(G'\cap \text{ker }\chi))$ is elementary Abelian for every $\chi\in \text{Irr}_{1}(G)$.  Since $[Z(\chi), G]\subseteq G'\cap \text{ker }\chi$, we get that $Z(\chi)/(G'\cap \text{ker }\chi)\subseteq Z(G/(G'\cap \text{ker }\chi))$. Consequently $Z(\chi)/(G'\cap \text{ker }\chi)$ is elementary Abelian and hence $Z(\chi)/ \text{ker }\chi$ is elementary Abelian. Since $Z(\chi)/ \text{ker }\chi$ is cyclic and $ \text{ker }\chi< Z(\chi)$, we get that $Z(\chi)/ \text{ker }\chi\cong C_{p}$. Since for every $\chi\in \text{Irr}_{1}(G)$, $\chi(1)^{2}=|G/Z(\chi)|=\frac{|G|}{p^{m+1}}$, we get that $|Z(\chi)|=p^{m+1}$. Thus we have that $|\text{ker }\chi|=p^{m}$ and hence kernel of each of the non-linear irreducible characters of $G$ has order $p^{m}$.
 \end{proof}
 \end{thm}


\end{document}